\documentclass[a4paper,11pt]{amsart}

\usepackage{graphicx}
\usepackage{mathptmx}
\usepackage{amsmath}
\usepackage{amssymb}
\usepackage{enumitem}
\usepackage{xcolor}

\newmuskip\pFqmuskip

\newcommand*\pFq[6][8]{%
  \begingroup 
  \pFqmuskip=#1mu\relax
  \mathcode`=\string"8000
  \begingroup\lccode`\~=`\,
  \lowercase{\endgroup\let~}\pFqcomma
  F^{#2}_{#3}{\left(\genfrac..{0pt}{}{#4}{#5}\bigg|#6\right)}%
  \endgroup
}
\newcommand{\pFqcomma}{\mskip\pFqmuskip}

\newtheorem{theorem}{Theorem}[section]

\newtheorem{remark}[theorem]{Remark}

\begin{document}

\title[]{Probabilistic degenerate poly-Bell polynomials associated with random variables}

\author{Pengxiang Xue}
\address{School of Science, Xi’an Technological University, Xi’an, 710021, Shaanxi, China}
\email{xuepx@xatu.edu.cn}
\author{Yuankui Ma*}
\address{School of Science, Xi’an Technological University, Xi’an, 710021, Shaanxi, China}
\email{mayuankui@xatu.edu.cn}
\author{Taekyun  Kim*}
\address{Department of Mathematics, Kwangwoon University, Seoul 139-701, Republic of Korea;School of Science, Xi’an Technological University, Xi’an, 710021, Shaanxi, China}
\email{tkkim@kw.ac.kr}
\author{Dae San  Kim*}
\address{Department of Mathematics, Sogang University, Seoul 121-742, Republic of Korea}
\email{dskim@sogang.ac.kr}
\author{Wenpeng Zhang*}
\address{School of Science, Xi’an Technological University, Xi’an, 710021, Shaanxi, China}
\email{wpzhang888@163.com}

\thanks{* corresponding authors}

\begin{abstract}
Let $Y$ be a random variable whose moment generating function exists in a neighborhood of the origin. The aim of this paper is to study the probabilistic degenerate poly-Bell polynomials associated with the random variable $Y$, arising from the degenerate polyexponential functions, which are probabilistic extensions of degenerate versions of the poly-Bell polynomials. We derive several explicit expressions and some related identities for them. In addition, we consider the special cases that $Y$ is the Bernoulli random variable with probability of success $p$ or the gamma random variable with parameters 1,1.
\end{abstract}

\maketitle

\markboth{\centerline{\scriptsize Probabilistic degenerate poly-Bell polynomials associated with random variables}}{\centerline{\scriptsize P. Xue,  Y. Ma, T. Kim,  D. S. Kim, W. Zhang}}

\section{Introduction}
Special functions and polynomials are indispensable tools in diverse fields such as mathematics, physics, engineering and related disciplines, addressing problems in areas like mathematical physics, numerical analysis, differential equations and quantum mechanics. Among these, probabilistic special polynomials, including probabilistic Bell, Fubini, and Stirling polynomials and numbers, have drawn significant attention recently. In [2], a probabilistic generalization of the Stirling numbers of the second kind was introduced, along with some of their properties. Subsequently, researchers have studied and explored numerous relations and properties of various probabilistic extensions of special polynomials, including Stirling polynomials of the second kind, degenerate Stirling polynomials of the second kind, Bernoulli polynomials, degenerate Bell polynomials, Lah numbers and polynomials and Bernstein polynomials (see [2,9,12,13,15,17,22-24,30-34] and the references therein). It is also remarkable that various degenerate versions of many special polynomials and numbers, gamma functions and umbral calculus have seen much attention with renewed interests and a lot of interesting results have been found (see [7,8,10,14,16,18,20,21,25-27] and the references therein). \par
In this paper, we investigate the probabilistic degenerate poly-Bell polynomials associated with $Y$, $\mathrm{Bel}_{n,\lambda}^{(k,Y)}(x)$, and find several explicit expressions and some related identities for them. We also consider the special cases that $Y$ is the Bernoulli random variable with probability of success $p$ or the gamma random variable with parameters 1,1. Here $Y$ is a random variable satisfying the moment condition in \eqref{1} below. The degenerate poly-Bell polynomials are defined by means of the degenerate polyexponential functions (see [11]) which are degenerate versions of the polyexponential functions (see [6]). Then the polynomials $\mathrm{Bel}_{n,\lambda}^{(k,Y)}(x)$ of our concern are probabilistic extensions of the degenerate poly-Bell polynomials (see \eqref{14}). \par
Assume that $Y$ is a random variable such that the moment generating function of $Y$,
\begin{equation}
E\big[e^{Yt}\big]=\sum_{n=0}^{\infty}E\big[Y^{n}\big]\frac{t^{n}}{n!},\quad (|t|<r), \label{1}
\end{equation}
exists for some $r>0$. Let $(Y_{j})_{j\ge 1}$ be a sequence of mutually independent copies of $Y$, and let
\begin{equation}
S_{k}=Y_{1}+Y_{2}+\cdots+Y_{k},\ (k\ge 1),\,\, \mathrm{with}\ S_{0}=0,\quad (\mathrm{see}\ [2]). \label{2}
\end{equation}

The outline of this paper is as follows. In Section 1, we recall the Stirling numbers of both kinds and the Bell polynomials. We remind the reader of the probabilistic Stirling numbers of the second associated with $Y$ and the probabilistic Bell polynomials associated with $Y$. We recall the degenerate exponentials, the degenerate Stirling numbers of both kinds, the degenerate logarithm and the degenerate Bell polynomials. We remind the reader of the probabilistic degenerate Stirling numbers of the second kind associated with $Y$, ${n \brace k}_{Y, \lambda}$, and the probabilistic Bell polynomials associated with $Y$. We recall the degenerate polyexponential  functions and the degenerate poly-Bell polynomials. We remind the reader of the gamma random variable with parameters $\alpha, \beta$, and the Lah numbers.
Section 2 contains the main results of this paper. We define the probabilistic degenerate poly-Bell polynomials $\mathrm{Bel}_{n,\lambda}^{(k,Y)}(x)$ (see \eqref{17}), arising from the degenerate polyexponential functions (see \eqref{13}). Here the polyexponential functions were first studied by Hardy [5] and reconsidered in [8], and their degenerate versions, namely the degenerate polyexponential functions, were introduced in [11]. We derive an explicit expression for $\mathrm{Bel}_{n,\lambda}^{(k,Y)}(x)$ as a finite sum involving ${n \brace l}_{Y,\lambda}$ in Theorem 2.2, and as a finite sum involving $E[(S_{m})_{n,\lambda}]$ (see \eqref{2}, \eqref{8}) in Theorem 2.3. We obtain identities involving $\mathrm{Bel}_{n,\lambda}^{(k,Y)}(x)$ in Theorems 2.4 and 2.5. In Theorems 2.6 and 2.7, we find explicit expressions for $\mathrm{Bel}_{n,\lambda}^{(k,Y)}(x)$ as finite sums involving the degenerate Stirling numbers of the second when $Y$ is the Bernoulli random variable with probability of success $p$ (see \eqref{4}). Finally, we get an explicit expression for $\mathrm{Bel}_{n,\lambda}^{(k,Y)}(x)$ when $Y$ is the gamma random variable with parameters 1, 1 (see \eqref{15}). \par
As general references, we let the reader refer to [1,3,4,28,29,33]. For the rest of this section, we recall what are needed throughout this paper. \par

\vspace{0.1in}

The Stirling numbers of the first kind are defined by
\begin{align}
&(x)_{n}=\sum_{k=0}^{n}S_{1}(n,k)x^{k},\quad (n\ge 0), \label{3} \\
&\frac{1}{k!}\big(\log(1+t)\big)^{k}=\sum_{n=k}^{\infty}S_{1}(n,k)\frac{t^{n}}{n!},
\quad (k \ge 0), \quad (\mathrm{see}\ [4,28]), \nonumber
\end{align}
where $(x)_{0}=1,\ (x)_{n}=x(x-1)\cdots(x-n+1),\ (n\ge 1)$. \\
As the inversion formula of \eqref{3}, the Stirling numbers of the second kind are defined by
\begin{align}
&x^{n}=\sum_{k=0}^{n}{n \brace k}(x)_{k},\quad (n\ge 0), \label{4} \\
&\frac{1}{k!}(e^{t}-1)^{k}=\sum_{n=k}^{\infty} {n \brace k}\frac{t^{n}}{n!}, \quad (k \ge 0), \quad (\mathrm{see}\ [4,28]). \nonumber
\end{align}
The Bell polynomials are defined by
\begin{equation}
\phi_{n}(x)=\sum_{k=0}^{n}{n \brace k}x^{k},\quad (n\ge 0),\quad (\mathrm{see}\ [1,4,11,28]). \label{5}
\end{equation}
When $x=1$, $\phi_{n}=\phi_{n}(1)$ are called the Bell numbers. \par
Recently, Adell introduced the probabilistic Stirling numbers of the second kind associated with $Y$, which are defined by
\begin{equation}
\frac{1}{k!}\Big(E\big[e^{Yt}\big]-1\Big)^{k}=\sum_{n=k}^{\infty}{n \brace k}_{Y}\frac{t^{n}}{n!},\quad (k\ge 0),\quad (\mathrm{see}\ [2]). \label{6}
\end{equation}
From \eqref{6}, we have
\begin{equation*}
{n \brace k}_{Y}=\frac{1}{k!}\sum_{j=0}^{k}\binom{k}{j}(-1)^{k-j}E\big[S_{j}^{n}\big],\ (0\le k\le n),\quad (\mathrm{see}\ [2]). 	
\end{equation*}
In [30], as a probabilisitc extension of the Bell polynomials (see \eqref{5}), the probabilistic Bell polynomials associated with $Y$ are given by
\begin{equation}
\sum_{n=0}^{\infty}\phi_{n}^{Y}(x)\frac{t^{n}}{n!}=e^{x( E[e^{Yt}]-1)}. \label{7}
\end{equation} \par
Let $\lambda$ be any nonzero real number. The degenerate falling factorial sequence is defined by
\begin{equation}
(x)_{0,\lambda}=1,\quad (x)_{n,\lambda}=x(x-\lambda)(x-2\lambda)\cdots\big(x-(n-1)\lambda\big),\ (n\ge 1). \label{8}
\end{equation}
The degenerate Stirling numbers of the second kind are defined by
\begin{align}
&(x)_{n,\lambda}=\sum_{k=0}^{n}{n \brace k}_{\lambda}(x)_{k},\quad (n\ge 0),\label{9}\\
&\frac{1}{k!}(e_{\lambda}(t)-1)^{k}=\sum_{n=k}^{\infty} {n \brace k}_{\lambda}\frac{t^{n}}{n!}, \quad (k \ge 0),\quad (\mathrm{see}\ [10,14,18,19]). \nonumber
\end{align}
As the inversion formula of \eqref{9}, the degenerate Stirling numbers of the first kind are given by
\begin{align}
&(x)_{n}=\sum_{k=0}^{n}S_{1,\lambda}(n,k)(x)_{k,\lambda},\quad (n\ge 0),\label{10} \\
&\frac{1}{k!}\big(\log_{\lambda}(1+t)\big)^{k}=\sum_{n=k}^{\infty}S_{1,\lambda}(n,k)\frac{t^{n}}{n!}, \quad (\mathrm{see}\ [10,14,18,19]). \nonumber
\end{align}
Note that
\begin{displaymath}
\lim_{\lambda\rightarrow 0}{n \brace k}_{\lambda}={n \brace k},\quad \lim_{\lambda\rightarrow 0}S_{1,\lambda}(n,k)=S_{1}(n,k).
\end{displaymath} \\
The degenerate Bell polynomials are given by
\begin{equation}
e^{x(e_{\lambda}(t)-1)}=\sum_{n=0}^{\infty}\phi_{n,\lambda}(x)\frac{t^{n}}{n!},\quad (\mathrm{see}\ [9,11]). \label{10-1}
\end{equation}
Note that
\begin{displaymath}
\phi_{n,\lambda}(x)=\sum_{k=0}^{n}{n \brace k}_{\lambda}x^{k},\quad (n \ge 0).
\end{displaymath} \par
Here the degenerate exponentials, $e_{\lambda}^{x}(t)$, are given by
\begin{equation}
e_{\lambda}^{x}(t)=\sum_{n=0}^{\infty}(x)_{n,\lambda}\frac{t^{n}}{n!}=(1+\lambda t)^{\frac{x}{\lambda}},\quad e_{\lambda}(t)=e_{\lambda}^{1}(t),\quad (\mathrm{see}\ [9,10,14,18,19]),  \label{11}
\end{equation}
and the degenerate logarithms, $\log_{\lambda}(t)$, are defined by
\begin{equation*}
\log_{\lambda}(1+t)=\sum_{n=1}^{\infty}\lambda^{n-1}(1)_{n,1/\lambda}\frac{t^{n}}{n!}=\frac{1}{\lambda}\big((1+t)^{\lambda}-1\big),\quad (\mathrm{see}\ [9,10,14,18,19]).
\end{equation*}
They are compositional inverses to each other so that
\begin{equation*}
e_{\lambda}\big(\log_{\lambda}(1+t)\big)=\log_{\lambda}\big(e_{\lambda}(1+t)\big)=1+t.
\end{equation*} \par
In [9], the probabilistic degenerate Stirling numbers of the second kind associated with $Y$, \,${n \brace k}_{Y,\lambda}$,\, and the probabilistic degenerate Bell polynomials associated with $Y$,\, $\phi_{n,\lambda}^{Y}(x)$,\, are respectively given by
\begin{equation}
\frac{1}{k!}\Big(E\big[e_{\lambda}^{Y}(t)\big]-1\Big)^{k}=\sum_{n=k}^{\infty}{n \brace k}_{Y,\lambda}\frac{t^{n}}{n!},\quad (k\ge 0), \label{12}
\end{equation}
and
\begin{equation*}
e^{x(E[e_{\lambda}^{Y}(t)]-1)}=\sum_{n=0}^{\infty}\phi_{n,\lambda}^{Y}(x)\frac{t^{n}}{n!}.
\end{equation*} \par
Recently, Kim-Kim introduced the degenerate polyexponential function given by
\begin{equation}
\mathrm{Ei}_{k,\lambda}(x)=\sum_{n=1}^{\infty}\frac{(1)_{n,\lambda}}{(n-1)!n^{k}}x^{n},\quad (k\in\mathbb{Z}),\quad (\mathrm{see}\ [19]). \label{13}
\end{equation}
Note taht
\begin{displaymath}
\mathrm{Ei}_{1,\lambda}(x)=e_{\lambda}(x)-1.
\end{displaymath}
In [19], the degenerate poly-Bell polynomials, $\mathrm{Bel}_{n,\lambda}^{(k)}(x)$, are defined by
\begin{equation}
\mathrm{Ei}_{k,\lambda}\Big(x\big(e_{\lambda}(t)-1\big)\Big)=\sum_{n=1}^{\infty}\mathrm{Bel}_{n,\lambda}^{(k)}(x)\frac{t^{n}}{n!},\quad \mathrm{Bel}_{0,\lambda}^{(k)}(x)=1. \label{14}
\end{equation} \par
A continuous random variable $Y$ whose density function is given by
\begin{equation}
f(y)=\left\{\begin{array}{ccc}
\frac{\beta e^{-\beta y}(\beta y)^{\alpha-1}}{\Gamma(\alpha)}, & \textrm{if $y\ge 0$},\\
0, & \textrm{if $y < 0$},
\end{array}\right.\quad (\mathrm{see}\ [12,29]), \label{15}
\end{equation}
for $\alpha,\beta>0$, is said to be the gamma random variable with parameters $\alpha,\beta$, which is denoted by $Y\sim\Gamma(\alpha,\beta)$. \par
The Lah numbers are given by
\begin{equation}
\frac{1}{k!}\bigg(\frac{t}{1-t}\bigg)^{k}=\sum_{n=k}^{\infty}L(n,k)\frac{t^{n}}{n!},\quad (k \ge 0),\quad (\mathrm{see}\ [4,8,20,28]).\label{16}
\end{equation}
Thus, by \eqref{16}, we get
\begin{equation*}
L(n,k)=\frac{n!}{k!}\binom{n-1}{k-1},\quad (n \ge k\ge 0).
\end{equation*}

\section{Probabilistic degenerate poly-Bell polynomials associated with random variables}
For $k\in\mathbb{Z}$, we define the {\it{probabilistic degenerate poly-Bell polynomials associated with $Y$}}, $\mathrm{Bel}_{n,\lambda}^{(k,Y)}(x)$, by
\begin{equation}
\mathrm{Ei}_{k,\lambda}\Big(x\big(E\big[e_{\lambda}^{Y}(t)\big]-1\big)\Big)=\sum_{n=1}^{\infty}\mathrm{Bel}_{n,\lambda}^{(k,Y)}(x)\frac{t^{n}}{n!},\quad \mathrm{Bel}_{0,\lambda}^{(k,Y)}(x)=1.\label{17}
\end{equation}
Here $\mathrm{Bel}_{n,\lambda}^{(k,Y)}= \mathrm{Bel}_{n,\lambda}^{(k,Y)}(1)$ are called the {\it{probabilistic degenerate poly-Bell numbers associated with $Y$}}. In particular, for $Y=1$, we have $\mathrm{Bel}_{n,\lambda}^{(k,1)}(x)= \mathrm{Bel}_{n,\lambda}^{(k)}(x)$. \par
For $k=1$, by using  \eqref{12} and \eqref{17}, we have
\begin{align}
\mathrm{Ei}_{1,\lambda}\Big(x\big(E\big[e_{\lambda}^{Y}(t)\big]-1\big)\Big)&=\sum_{n=1}^{\infty}(1)_{n,\lambda}x^{n}\frac{1}{n!}\Big(E\big[e_{\lambda}^{Y}(t)\big)-1\Big)^{n} \label{18}\\
&=\sum_{n=1}^{\infty}(1)_{n,\lambda}x^{n}\sum_{m=n}^{\infty}{m \brace n}_{Y,\lambda}\frac{t^{m}}{m!}\nonumber\\
&=\sum_{m=1}^{\infty}\sum_{n=1}^{m}(1)_{n,\lambda}{m \brace n}_{Y,\lambda}x^{n}\frac{t^{m}}{m!}. \nonumber
\end{align}
Therefore, by \eqref{17} and \eqref{18}, we obtain the following theorem.
\begin{theorem}
For $m\in\mathbb{Z}$, we have
\begin{displaymath}
\mathrm{Bel}_{m,\lambda}^{(1,Y)}(x)=\sum_{n=1}^{m}(1)_{n,\lambda}{m \brace n}_{Y,\lambda}x^{n}.
\end{displaymath}
In particular, letting $x=1$ gives
\begin{displaymath}
\mathrm{Bel}_{n,\lambda}^{(1,Y)}=\sum_{n=1}^{m}(1)_{n,\lambda}{m \brace n}_{Y,\lambda}.
\end{displaymath}
\end{theorem}
From \eqref{12} and \eqref{17}, we note that
\begin{align}
\sum_{n=1}^{\infty} \mathrm{Bel}_{n,\lambda}^{(k,Y)}(x)\frac{t^{n}}{n!}&= \mathrm{Ei}_{k,\lambda}\Big(x\big(E\big[e_{\lambda}^{Y}(t)\big]-1\big)\Big)\label{19}\\
&=\sum_{l=1}^{\infty}\frac{(1)_{l,\lambda}x^{l}}{l^{k-1}}\frac{1}{l!}\Big(E\big[e_{\lambda}^{Y}(t)\big]-1\Big)^{l}\nonumber\\
&=\sum_{l=1}^{\infty}\frac{(1)_{l,\lambda}x^{l}}{l^{k-1}}\sum_{n=l}^{\infty}{n \brace l}_{Y,\lambda}\frac{t^{n}}{n!}\nonumber\\
&=\sum_{n=1}^{\infty}\sum_{l=1}^{n}\frac{(1)_{l,\lambda}}{l^{k-1}}{n \brace l}_{Y,\lambda}x^{l}\frac{t^{n}}{n!}. \nonumber
\end{align}
Therefore, by \eqref{19}, we obtain the following theorem.
\begin{theorem}
For $n\in\mathbb{N}$, we have
\begin{displaymath}
\mathrm{Bel}_{n,\lambda}^{(k,Y)}(x)=\sum_{l=1}^{n}\frac{(1)_{l,\lambda}}{l^{k-1}}{n \brace l}_{Y,\lambda}x^{l}.
\end{displaymath}
\end{theorem}
By \eqref{2}, \eqref{8} and \eqref{17}, we get
\begin{align}
\sum_{n=1}^{\infty}\mathrm{Bel}_{n,\lambda}^{(k,Y)}(x)\frac{t^{n}}{n!}&=\sum_{l=1}^{\infty}\frac{(1)_{l,\lambda}x^{l}}{(l-1)!l^{k}}\Big(E\big[e_{\lambda}^{Y}(t)\big]-1\Big)^{l}\label{20}\\
&=\sum_{l=1}^{\infty}\frac{(1)_{l,\lambda}x^{l}}{(l-1)!l^{k}}\sum_{m=0}^{l}\binom{l}{m}(-1)^{l-m}\Big(E\big[e_{\lambda}^{Y}(t)\big]\Big)^{m}\nonumber\\
&=\sum_{l=1}^{\infty}\sum_{m=0}^{l}\frac{(1)_{l,\lambda}x^{l}	}{(l-1)!l^{k}}\binom{l}{m}(-1)^{l-m}E\Big[e_{\lambda}^{Y_{1}+Y_{2}+\cdots+Y_{m}}(t)\Big] \nonumber \\
&=\sum_{n=0}^{\infty}\sum_{l=1}^{\infty}\sum_{m=0}^{l}\frac{(1)_{l,\lambda}x^{l}}{(l-1)!l^{k}}\binom{l}{m}(-1)^{l-m}E\big[(S_{m})_{n,\lambda}\big]\frac{t^{n}}{n!}\nonumber\\
&=\sum_{n=1}^{\infty}\sum_{l=1}^{\infty}\sum_{m=0}^{l}\frac{(1)_{l,\lambda}x^{l}}{(l-1)!l^{k}}\binom{l}{m}(-1)^{l-m}E\big[(S_{m})_{n,\lambda}\big]\frac{t^{n}}{n!}.\nonumber
\end{align}
Therefore, by \eqref{20} and noting from Theorem 2.2 that $\mathrm{Bel}_{n,\lambda}^{(k,Y)}(x)$ is a polynomial of degree $n$, we obtain the following theorem.
\begin{theorem}
For $n\in\mathbb{N}$, we have
\begin{equation*}
\mathrm{Bel}_{n,\lambda}^{(k,Y)}(x)= \sum_{l=1}^{n}\bigg(\frac{(1)_{l,\lambda}}{(l-1)!l^{k}}\sum_{m=0}^{l}\binom{l}{m}(-1)^{l-m}E\big[(S_{m})_{n,\lambda}\big]\bigg)x^{l},
\end{equation*}
and
\begin{equation*}
(1)_{l,\lambda}\sum_{m=0}^{l}\binom{l}{m}(-1)^{l-m}E\big[(S_{m})_{n,\lambda}\big]=0, \quad \mathrm{for}\,\,\, l \ge n+1.
\end{equation*}
\end{theorem}
\begin{remark}
From Theorem 2.3, we have
\begin{equation}
\sum_{m=0}^{l}\binom{l}{m}(-1)^{l-m}E\big[(S_{m})_{n,\lambda}\big]=0, \label{20-1}
\end{equation}
for\,\, $l \ge n+1 \ge 2$\,\,and\,\, $\lambda \ne 1, \frac{1}{2}, \dots, \frac{1}{l-1}$. \\
Furthermore, letting $\lambda \rightarrow 0$ in \eqref{20-1} yields
\begin{equation}
\sum_{m=0}^{l}\binom{l}{m}(-1)^{l-m}E\big[S_{m}^{n}\big]=0, \quad \mathrm{for}\,\,\, l \ge n+1 \ge 2. \label{20-2}
\end{equation}
Let $Y$ be the gamma random variable with parameters $\alpha, \beta >0$ (see \eqref{15}). Then we see that
\begin{equation*}
E[e^{Y t}]=\Big(\frac{\beta}{\beta-t} \Big)^{\alpha} ,\quad (\mathrm{see}\ [29], \ p.62).
\end{equation*}
Note that
\begin{equation}
\sum_{n=0}^{\infty}E[S_{m}^{n}]\frac{t^{n}}{n!}=E[e^{Y t}]^{m}=\Big(\frac{\beta}{\beta-t} \Big)^{\alpha m}=\sum_{n=0}^{\infty}\frac{(\alpha m+n-1)_{n}}{\beta^{n}}\frac{t^{n}}{n!}. \label{20-3}
\end{equation}
From \eqref{20-2} and \eqref{20-3}, we obtain the following identity:
\begin{equation*}
\sum_{m=0}^{l}\binom{l}{m}(-1)^{l-m}(\alpha m +n-1)_{n}=0, \quad \mathrm{for}\,\,\, l \ge n+1 \ge 2.
\end{equation*}
One could find many interesting identities by using \eqref{20-1} and \eqref{20-2} and making different choices for $Y$.
\end{remark}

Replacing $t$ by $\log_{\lambda}(1+t)$ in \eqref{17} and using \eqref{10}, we have
\begin{align}
\mathrm{Ei}_{k,\lambda}\Big(x\big(E[(1+t)^{Y}]-1\big)\Big)&=\sum_{l=1}^{\infty}\mathrm{Bel}_{l,\lambda}^{(k,Y)}(x)\frac{1}{l!}\Big(\log_{\lambda}(1+t)\Big)^{l}\label{21}\\
&=\sum_{l=1}^{\infty}\mathrm{Bel}_{l,\lambda}^{(k,Y)}(x)\sum_{n=l}^{\infty}S_{1,\lambda}(n,l)\frac{t^{n}}{n!} \nonumber \\
&=\sum_{n=1}^{\infty}\sum_{l=1}^{n}\mathrm{Bel}_{l,\lambda}^{(k,Y)}(x)S_{1,\lambda}(n,l)\frac{t^{n}}{n!}.\nonumber	
\end{align}
On the other hand, by utilizing \eqref{3}, \eqref{6} and \eqref{17}, we also have
\begin{align}
\mathrm{Ei}_{k,\lambda} \Big(x\big(E[(1+t)^{Y}]-1\big)\Big)&=\sum_{l=1}^{\infty}\frac{(1)_{l,\lambda}x^{l}}{l^{k}(l-1)!}\Big(E\big[(1+t)^{Y}\big]-1\Big)^{l}\label{22}\\
&=\sum_{l=1}^{\infty}\frac{x^{l}(1)_{l,\lambda}}{l^{k-1}}\frac{1}{l!}\Big(E\big[e^{Y\log (1+t)}\big]-1\Big)^{l} \nonumber \\
&=\sum_{l=1}^{\infty}\frac{x^{l}(1)_{l,\lambda}}{l^{k-1}}\sum_{j=l}^{\infty}{j \brace l}_{Y}\frac{1}{j!}\big(\log(1+t)\big)^{j}\nonumber\\
&=\sum_{j=1}^{\infty}\sum_{l=1}^{j}\frac{(1)_{l,\lambda}x^{l}}{l^{k-1}}{j \brace l}_{Y}\sum_{n=j}^{\infty}S_{1}(n,l)\frac{t^{n}}{n!}\nonumber\\
&=\sum_{n=1}^{\infty}\sum_{j=1}^{n}\sum_{l=1}^{j}\frac{(1)_{l,\lambda}}{l^{k-1}}{j \brace l}_{Y}S_{1}(n,l)x^{l}\frac{t^{n}}{n!}.\nonumber
\end{align}
Therefore, by \eqref{21} and \eqref{22}, we obtain the following theorem.
\begin{theorem}
For $n\in\mathbb{N}$, we have
\begin{displaymath}
\sum_{l=1}^{n}S_{1,\lambda}(n,l)\mathrm{Bel}_{l,\lambda}^{(k,Y)}(x)=\sum_{j=1}^{n}\sum_{l=1}^{j}\frac{(1)_{l,\lambda}}{l^{k-1}}{j \brace l}_{Y}S_{1}(n,l)x^{l}.
\end{displaymath}
\end{theorem}
From \eqref{13} and \eqref{17}, we note that
\begin{align}
&\frac{\partial}{\partial t} \mathrm{Ei}_{k,\lambda} \Big(x\big(E[e_{\lambda}^{Y}(t)]-1\big)\Big)=\frac{\partial}{\partial t}\sum_{n=1}^{\infty}\frac{(1)_{n,\lambda}x^{n}\big(E[e_{\lambda}^{Y}(t)]-1\big)^{n}}{(n-1)!n^{k}}\label{23}\\
&=\sum_{n=1}^{\infty}\frac{(1)_{n,\lambda}x^{n}}{(n-1)!n^{k-1}}\Big(E\big[e_{\lambda}^{Y}(t)\big]-1\Big)^{n-1}E\big[Ye_{\lambda}^{Y-\lambda}(t)\big], \nonumber
\end{align}
and
\begin{align}
\frac{\partial}{\partial t} \mathrm{Ei}_{k,\lambda} \Big(x\big(E[e_{\lambda}^{Y}(t)]-1\big)\Big)&=\frac{\partial}{\partial t}\sum_{n=1}^{\infty}\mathrm{Bel}_{n,\lambda}^{(k,Y)}(x)\frac{t^{n}}{n!} \label{24}\\
&=\sum_{n=0}^{\infty}\mathrm{Bel}_{n+1,\lambda}^{(k,Y)}(x)\frac{t^{n}}{n!}.\nonumber
\end{align}
By \eqref{11}, \eqref{13}, \eqref{17}, \eqref{23} and \eqref{24}, we get
\begin{align}
\big(E\big[e_{\lambda}^{Y}(t)\big]-1\big)\sum_{m=0}^{\infty}\mathrm{Bel}_{m+1,\lambda}^{(k,Y)}(x)\frac{t^{m}}{m!}&=\sum_{n=1}^{\infty}\frac{x^{n}(1)_{n,\lambda}}{(n-1)!n^{k-1}}\big(E\big[e_{\lambda}^{Y}(t)\big]-1\big)^{n}E\big[Ye_{\lambda}^{Y-\lambda}(t)\big]\label{25}\\
&=\mathrm{Ei}_{k-1,\lambda}\big(x\big(E\big[e_{\lambda}^{Y}(t)\big]-1\big)\big)E\big[Ye_{\lambda}^{Y-\lambda}(t)\big] \nonumber	\\
&=\sum_{m=1}^{\infty}\mathrm{Bel}_{m,\lambda}^{(k-1,Y)}(x)\frac{t^{m}}{m!}\sum_{l=0}^{\infty}E\big[Y(Y-\lambda)_{l,\lambda}\big]\frac{t^{l}}{l!}\nonumber\\
&=\sum_{n=1}^{\infty}\sum_{m=1}^{n}\binom{n}{m}\mathrm{Bel}_{m,\lambda}^{(k-1,Y)}(x)E\big[Y(Y-\lambda)_{n-m,\lambda}\big]\frac{t^{n}}{n!}.\nonumber
\end{align}
On the other hand, by using \eqref{11}, we also have
\begin{align}
\big(E\big[e_{\lambda}^{Y}(t)\big]-1\big)\sum_{m=0}^{\infty}\mathrm{Bel}_{m+1,\lambda}^{(k,Y)}(x)\frac{t^{m}}{m!}&=\sum_{l=1}^{\infty}E\big[(Y)_{l,\lambda}\big]\frac{t^{l}}{l!}\sum_{m=0}^{\infty}\mathrm{Bel}_{m+1,\lambda}^{(k,Y)}(x)\frac{t^{m}}{m!}\label{26}\\
&=\sum_{n=1}^{\infty}\sum_{m=0}^{n-1}\binom{n}{m}\mathrm{Bel}_{m+1,\lambda}^{(k,Y)}(x)E\big[(Y)_{n-m,\lambda}\big]\frac{t^{n}}{n!}.\nonumber
\end{align}
Therefore, by \eqref{25} and \eqref{26}, we obtain the following theorem.
\begin{theorem}
For $n\in\mathbb{N}$, we have
\begin{displaymath}
\sum_{m=1}^{n}\binom{n}{m}\mathrm{Bel}_{m,\lambda}^{(k-1,Y)}(x)E\big[Y(Y-\lambda)_{n-m,\lambda}\big]=\sum_{m=0}^{n-1}\binom{n}{m}\mathrm{Bel}_{m+1,\lambda}^{(k,Y)}(x)E\big[(Y)_{n-m,\lambda}\big].
\end{displaymath}
\end{theorem}
Let $Y$ be the Bernoulli random variable with probability of success $p$. Noting that $E[e_{\lambda}^{Y}(t)]=1+p\big(e_{\lambda}(t)-1\big)$, and using \eqref{9}, \eqref{13} and \eqref{17}, we have
\begin{align}
\sum_{n=1}^{\infty}\mathrm{Bel}_{n,\lambda}^{(k,Y)}(x)\frac{t^{n}}{n!}	&=\sum_{m=1}^{\infty}\frac{(1)_{m,\lambda}x^{m}}{m^{k}(m-1)!}\Big(E\big[e_{\lambda}^{Y}(t)\big]-1\Big)^{m} \label{27}\\
&=\sum_{m=1}^{\infty}\frac{(1)_{m,\lambda}x^{m}}{m^{k-1}}p^{m}\frac{1}{m!}\big(e_{\lambda}(t)-1\big)^{m}\nonumber\\
&=\sum_{m=1}^{\infty}\frac{(1)_{m,\lambda}x^{m}}{m^{k-1}}p^{m}\sum_{n=m}^{\infty}{n \brace m}_{\lambda}\frac{t^{n}}{n!} \nonumber \\
&=\sum_{n=1}^{\infty}\sum_{m=1}^{n}\frac{(1)_{m,\lambda}p^{m}}{m^{k-1}}{n \brace m}_{\lambda}x^{m}\frac{t^{n}}{n!}. \nonumber
\end{align}
Therefore, by \eqref{27}, we obtain the following theorem.
\begin{theorem}
Let $Y$ be the Bernoulli random variable with probability of success $p$. For $n\in\mathbb{N}$, we have
\begin{equation*}
\mathrm{Bel}_{n,\lambda}^{(k,Y)}(x)=\sum_{m=1}^{n}\frac{(1)_{m,\lambda}}{m^{k-1}}p^{m}{n \brace m}_{\lambda}x^{m}.
\end{equation*}
In particular, letting $x=1$ yields
\begin{equation*}
\mathrm{Bel}_{n,\lambda}^{(k,Y)}=\sum_{m=1}^{n}\frac{(1)_{m,\lambda}}{m^{k-1}}p^{m}{n \brace m}_{\lambda}.
\end{equation*}
\end{theorem}
From \eqref{10} and \eqref{13}, we note that
\begin{align}
\mathrm{Ei}_{k,\lambda}\big(\log_{\lambda}(1+t)\big)&=\sum_{m=1}^{\infty}\frac{(1)_{m,\lambda}}{m^{k-1}}\frac{1}{m!}\Big(\log_{\lambda}(1+t)\Big)^{m}\label{28}\\
&=\sum_{m=1}^{\infty}\frac{(1)_{m,\lambda}}{m^{k-1}}\sum_{n=m}^{\infty}S_{1,\lambda}(n,m)\frac{t^{n}}{n!}\nonumber\\
&=\sum_{n=1}^{\infty}\sum_{m=1}^{n}\frac{(1)_{m,\lambda}}{m^{k-1}}S_{1,\lambda}(n,m)\frac{t^{n}}{n!}.\nonumber
\end{align}
Replacing $t$ by $e_{\lambda}\Big(x\big(E\big[e_{\lambda}^{Y}(t)\big]-1\big)\Big)-1$ in \eqref{28} and utilizing \eqref{9} and \eqref{12}, we have
\begin{align}
\mathrm{Ei}_{k,\lambda}\Big(x\big(E\big[e_{\lambda}^{Y}(t)\big]-1\big)\Big)&=\sum_{l=1}^{\infty}\sum_{m=1}^{l}\frac{(1)_{m,\lambda}}{m^{k-1}}S_{1,\lambda}(l,m)\frac{1}{l!}\Big(e_{\lambda}\big(x\big(E[e_{\lambda}^{Y}(t)\big)-1\big)-1\Big)^{l} \label{29}\\
&=\sum_{l=1}^{\infty}\sum_{m=1}^{l}\frac{(1)_{m,\lambda}}{m^{k-1}}S_{1,\lambda}(l,m)\sum_{j=l}^{\infty}{j \brace l}_{\lambda}\frac{1}{j!}x^{j}\Big(E\big[e_{\lambda}^{Y}(t)\big]-1\Big)^{j}\nonumber\\
&=\sum_{j=1}^{\infty}\sum_{l=1}^{j}\sum_{m=1}^{l}\frac{(1)_{m,\lambda}}{m^{k-1}}S_{1,\lambda}(l,m){j \brace l}_{\lambda}x^{j}\frac{1}{j!}\Big(E\big[e_{\lambda}^{Y}(t)\big]-1\Big)^{j}.\nonumber
\end{align} \par
Assume that $Y$ is the Bernoulli random variable with probability of success $p$. Thus, by \eqref{9} and \eqref{29}, we get
\begin{align}
\sum_{n=1}^{\infty}\mathrm{Bel}_{n,\lambda}^{(k,Y)}(x)\frac{t^{n}}{n!}&=\sum_{j=1}^{\infty}\sum_{l=1}^{j}\sum_{m=1}^{l}\frac{(1)_{m,\lambda}}{m^{k-1}}S_{1,\lambda}(l,m){j \brace l}_{\lambda}p^{j}x^{j}\frac{1}{j!}\big(e_{\lambda}(t)-1\big)^{j}\label{30}\\
&=\sum_{j=1}^{\infty}\sum_{l=1}^{j}\sum_{m=1}^{l}\frac{(1)_{m,\lambda}}{m^{k-1}}S_{1,\lambda}(l,m){j\brace l}_{\lambda}p^{j}x^{j}\sum_{n=j}^{\infty}{n \brace j}_{\lambda}\frac{t^{n}}{n!}\nonumber\\
&=\sum_{n=1}^{\infty}\sum_{j=1}^{n}\sum_{l=1}^{j}\sum_{m=1}^{l}\frac{(1)_{m,\lambda}}{m^{k-1}}S_{1,\lambda}(l,m){j \brace l}_{\lambda}{n \brace j}_{\lambda}p^{j}x^{j}\frac{t^{n}}{n!}.\nonumber
\end{align}
Therefore, by \eqref{30}, we obtain the following theorem.
\begin{theorem}
Let $Y$ be the Bernoulli random variable with probability of success $p$. For $n\in\mathbb{N}$, we have
\begin{displaymath}
\mathrm{Bel}_{n,\lambda}^{(k,Y)}(x)= \sum_{j=1}^{n}\sum_{l=1}^{j}\sum_{m=1}^{l}\frac{(1)_{m,\lambda}}{m^{k-1}}S_{1,\lambda}(l,m){j \brace l}_{\lambda}{n \brace j}_{\lambda}p^{j}x^{j}.
\end{displaymath}
\end{theorem}
Let $Y\sim\Gamma(1,1)$,\,\, (see \eqref{15}). Then we see that
\begin{equation}
E[e_{\lambda}^{Y}(t)]=\int_{0}^{\infty}e_{\lambda}^{y}(t)e^{-y}dy=\frac{1}{1-\frac{1}{\lambda}\log(1+\lambda t)},\quad \Big(t < \frac{e^{\lambda}-1}{\lambda}\Big). \label{31}
\end{equation}
By using \eqref{3}, \eqref{13}, \eqref{16}, \eqref{17} and \eqref{31}, we proceed
\begin{align}
\sum_{n=1}^{\infty}\mathrm{Bel}_{n,\lambda}^{(k,Y)}(x)\frac{t^{n}}{n!}&=\sum_{m=1}^{\infty}\frac{(1)_{m,\lambda}x^{m}}{m^{k}(m-1)!}\Big(E\big[e_{\lambda}^{Y}(t)\big]-1\Big)^{m} \label{32}\\
&=\sum_{m=1}^{\infty}\frac{(1)_{m,\lambda}x^{m}}{m^{k-1}}\frac{1}{m!}\bigg(\frac{\frac{1}{\lambda}\log(1+\lambda t)}{1-\frac{1}{\lambda}\log(1+\lambda t) }\bigg)^{m}\nonumber\\
&=\sum_{m=1}^{\infty}\frac{(1)_{m,\lambda}x^{m}}{m^{k-1}}\sum_{l=m}^{\infty}L(l,m)\frac{1}{\lambda^{l}}\frac{1}{l!}\big(\log(1+\lambda t)\big)^{l} \nonumber \\
&=\sum_{m=1}^{\infty}\frac{(1)_{m,\lambda}x^{m}}{m^{k-1}}\sum_{l=m}^{\infty}L(l,m)\frac{1}{\lambda^{l}}\sum_{n=l}^{\infty}S_{1}(n,l)\lambda^{n}\frac{t^{n}}{n!}
\nonumber \\
&=\sum_{n=1}^{\infty}\sum_{m=1}^{n}\sum_{l=m}^{n}\frac{(1)_{m,\lambda}}{m^{k-1}}\lambda^{n-l}L(l,m)S_{1}(n,l)x^{m}\frac{t^{n}}{n!}. \nonumber
\end{align}
Therefore, by \eqref{32}, we obtain the following theorem.
\begin{theorem}
Let $Y\sim\Gamma(1,1)$. Then, for $n\in\mathbb{N}$, we have
\begin{displaymath}
\mathrm{Bel}_{n,\lambda}^{(k,Y)}(x)= \sum_{m=1}^{n}\sum_{l=m}^{n}\frac{(1)_{m,\lambda}}{m^{k-1}}\lambda^{n-l}L(l,m)S_{1}(n,l)x^{m}.
\end{displaymath}
In particular, letting $x=1$ gives
\begin{displaymath}
\mathrm{Bel}_{n,\lambda}^{(k,Y)}= \sum_{m=1}^{n}\sum_{l=m}^{n}\frac{(1)_{m,\lambda}}{m^{k-1}}\lambda^{n-l}L(l,m)S_{1}(n,l).
\end{displaymath}
\end{theorem}

\section{Conclusion}
We investigated the probabilistic degenerate poly-Bernoulli polynomials associated with $Y$, $\mathrm{Bel}_{n,\lambda}^{(k,Y)}(x)$, (see {\eqref{13}, \eqref{17}), and derived several explicit expressions and related identities for them.
Notably, as it was noted in Remark 2.4, we obtained the following interesting identity:
\begin{equation}
\sum_{m=0}^{l}\binom{l}{m}(-1)^{l-m}E\big[(S_{m})_{n,\lambda}\big]=0, \label{33}
\end{equation}
for\,\, $l \ge n+1 \ge 2$\,\,and\,\, $\lambda \ne 1, \frac{1}{2}, \dots, \frac{1}{l-1}$. \par
Let $Y$ be the Poisson random variable with parameter $\alpha > 0$. Then the probability mass function of $Y$ is given by
\begin{equation}
P\{Y=i \}=e^{-\alpha}\frac{\alpha^{i}}{i!},\quad i=0,1,2,\dots. \label{34}
\end{equation}
Then, by using \eqref{34}, one shows that
\begin{equation}
E[e_{\lambda}^{Y}(t)]=e^{\alpha \big(e_{\lambda}(t)-1\big)} ,\quad (\mathrm{see}\ [29], \ p.61). \label{35}
\end{equation}
Thus, by utilizing \eqref{35}, we have
\begin{equation}
\sum_{n=0}^{\infty}E[(S_{m})_{n,\lambda}]\frac{t^{n}}{n!}=E[e_{\lambda}^{Y}(t)]^{m}=e^{\alpha m\big(e_{\lambda}(t)-1 \big)}=\sum_{n=0}^{\infty}\phi_{n,\lambda}(\alpha m)\frac{t^{n}}{n!},\quad (\mathrm{see}\,\, \eqref{10-1}). \label{36}
\end{equation}
From \eqref{33} and \eqref{36}, we obtain the next identity:
\begin{equation*}
\sum_{m=0}^{l}\binom{l}{m}(-1)^{l-m}\phi_{n,\lambda}(\alpha m)=0,
\end{equation*}
for\,\, $l \ge n+1 \ge 2$\,\,and\,\, $\lambda \ne 1, \frac{1}{2}, \dots, \frac{1}{l-1}$.

\vspace{0.2in}

{\bf Acknowledgements}

The third author of this research was conducted by the Research Grant of Kwangwoon University in 2025.
This research was funded by the National Natural Science Foundation of China (No. 12271320).

\end{document}